\documentclass[11pt]{article}
\usepackage{amsmath,amsfonts,amssymb,latexsym}
\newtheorem{theorem}{Theorem}[section]
\newtheorem{proposition}[theorem]{Proposition}
\newtheorem{lemma}[theorem]{Lemma}
\newtheorem{cor}{Corollary}[theorem]
\def\R{{\mathbb R}}

\newcommand{\xR}{{]}{-\infty},+\infty]}
\newcommand{\Rex}{\xR}

\newcommand{\ps}{\smallbreak}

\newcommand{\proof}{\noindent{\bf Proof. }}
\newcommand{\cqfd}{\mbox{}\nolinebreak\hfill\rule{2mm}{2mm}\medbreak\par}

\newcommand{\lsc}{lower semicontinuous}

\newcommand{\del}{\partial}

\newcommand{\Sd}{\breve{\del}}
\newcommand{\dom} {{\rm dom} \kern.15em}
\newcommand{\tq}{:}
\newcommand{\la}{\langle}
\newcommand{\ra}{\rangle}
\newcommand{\ld}{\lambda}
\newcommand{\eps}{\varepsilon}

\newcommand{\bx}{\bar{x}}
\newcommand{\xb}{\bar{x}}
\newcommand{\xt}{x_\eps}

\begin{document}
\thispagestyle{empty}

\begin{center}
{\large\bf\sc Subdifferential Test for Optimality}
\end{center}

\begin{center}
  {\small\begin{tabular}{c}
  Florence JULES and Marc LASSONDE\\
   Universit\'e des Antilles et de la Guyane\\
  97159 Pointe \`a Pitre, France\\
  E-mail: florence.jules@univ-ag.fr, marc.lassonde@univ-ag.fr
  \end{tabular}}
\end{center}

\medbreak\noindent
\textbf{Abstract.}
We provide a first-order necessary and sufficient condition for
optimality of lower semicontinuous functions on Banach spaces
using the concept of subdifferential.
From the sufficient condition we derive that any subdifferential operator 
is monotone absorbing, hence maximal monotone when the function is convex.

\medbreak\noindent
\textbf{Keywords:}
  lower semicontinuity,
  subdifferential, directional derivative,
  first-order condition,
  optimality criterion, maximal monotonicity.
  
\medbreak\noindent
\textbf{2010 Mathematics  Subject Classification:}
  49J52, 49K27, 26D10, 26B25.
%%%%%%%%%%%%%%%%%%%%%%%%%%%%%%%%%%%%%%%%%%%%%%%%%%%%%%%%%%%%%%%%%%%%%%%
\section{Introduction}\label{intro}
First-order sufficient conditions for optimality in terms of derivatives
or directional derivatives are well known.
Typical such conditions are variants of Minty variational inequalities.
Let us quote for example the following simple result (see \cite{Bot86}):
{\it
Let $f$ be a real-valued function which is continuous on a neighbourhood $D$
centred at $a$ in $\R^n$ and differentiable on $D\setminus \{a\}$. Then,
$f$ has a local minimum value at $a$ if $(x-a)\cdot \nabla f(x)>0$
for all $x\in D\setminus \{a\}$.}
For first-order conditions in terms of directional derivatives, we refer to \cite{CGR04}.

The main objective of this note is to establish a necessary and sufficient condition
for optimality of nonsmooth \lsc\ functions on Banach spaces using subdifferentials.
To this end, we first provide 
a link between the directional derivative of a function and
its dual companion represented by the subdifferential (Theorem \ref{link}).
Then, we prove a sharp version of the mean value inequality
for \lsc\ function using directional derivative (Lemma \ref{mvi}). Finally,
we combine both results to obtain our subdifferential test for optimality
(Theorem \ref{main}).
A discussion on the sufficient condition follows where it is shown that
any subdifferential operator is \textit{monotone absorbing}, a property which reduces to 
maximal monotonicity when the function is convex
(Theorem \ref{mono_absorb}).

This paper complements our previous work \cite{JL12} which concerned only
the elementary subdifferentials.
%%%%%%%%%%%%%%%%%%%%%%%%%%%%%%%%%%%%%%%%%%%%%%%%%%%%%%%%%%
\section{Directional Derivative and Subdifferential}

In the following, $X$ is a real Banach space with unit ball $B_X$,
$X^*$ is its topological dual with unit ball $B_{X^*}$,
and $\la .,. \ra$ is the duality pairing.
Set-valued operators $T:X\rightrightarrows X^*$
are identified with their graph $T\subset X\times X^*$.
For a subset $A\subset X$, $x\in X$ and $\ld>0$, we let
$d_A(x):=\inf_{y\in A} \|x-y\|$ and $B_\ld(A):=\{ y\in X\tq d_A(y)\le \ld\}$. 
All the functions $f : X\to\Rex$ are assumed to be lower semicontinuous
and \textit{proper}, which means that
the set  $\dom f:=\{x\in X\tq f(x)<\infty\}$ is nonempty.
The (radial or lower Dini) \textit{directional derivative} of a function $f$ at a point
$\bx\in \dom f$ is given by:
\begin{equation*}\label{general-dd}
\forall d\in X,\quad f'(\xb;d):=\liminf_{t\searrow 0}\,\frac{f(\xb+td)-f(\xb)}{t}.
\end{equation*}
A \textit{subdifferential} of a \lsc\ function $f$ is a set-valued operator
$\partial f: X \rightrightarrows X^\ast$ which
coincides with the usual convex subdifferential whenever $f$ is convex, that is,
\[
\del f(\xb):= \{ x^* \in X^* \tq \la x^*,y-\xb\ra + f(\xb) \leq f(y),\, \forall y \in X \},\]
and satisfies elementary stability properties like
$\del(f-x^*)(x)=\del f(x) -x^*$ for every $x\in X$ and $x^*\in X^*$.
In this work, we also require that the subdifferentials satisfy
the following basic calculus rule: 
\medbreak
\textit{Separation Principle}.
For any \lsc\ $f,\varphi$ with $\varphi$ convex Lipschitz
near $\xb\in\dom f \cap\dom \varphi$,
if $f+\varphi$ admits a local minimum at $\xb$, then
$0\in \del f(\xb)+ \del \varphi(\xb).$
\ps
\textit{Examples}.
The Clarke subdifferential \cite{Cla90}, the Michel-Penot subdifferential \cite{MP92}, %Roc79,
the Ioffe A-subdifferential \cite{Iof89} satisfy the Separation Principle in any Banach space.
The limiting versions of the elementary subdifferentials
(proximal, Fr\'echet, Hadamard, G\^ateaux, \ldots) satisfy the Separation Principle
in appropriate Banach spaces (see, e.g., \cite{JL12,BZ05,Mor06ab} and the references therein).
\ps
The link between the directional derivative and the subdifferential
is described via the following $\eps$-enlargement of the subdifferential:
\begin{equation*}\label{FJ-sdiff}
\Sd_\eps f(\xb):=\{x^*_\eps\in X^*\tq x^*_\eps\in\del f(x_\eps) \mbox{ with }
\|x_\eps-\xb\|\le\eps,\,|f(x_\eps)-f(\xb)|\le\eps,\, \la x^*_\eps,x_\eps-\xb\ra\le \eps\}.
\end{equation*}

\begin{theorem}\label{link}
Let $X$ be a Banach space, $f:X\to\xR$ be lower semicontinuous, 
$\xb\in\dom f$ and $d\in X$. Then, for every $\eps>0$,
the sets $\Sd_\eps f(\xb)$ are nonempty and
\begin{equation}\label{CDD-formula}
%f'(\xb;d)\le\inf_{\eps>0}\sup\{\la x^*_\eps,d\ra\tq x^*_\eps\in \Sd_\eps f(\xb)\}.
f'(\xb;d)\le\inf_{\eps>0}\sup\la \Sd_\eps f(\xb),d\ra.
\end{equation}
\end{theorem}
\proof
We first show that the sets $\Sd_\eps f(\xb)$ are nonempty.
The arguments are standard, we give them for completeness.
Since $f$ is \lsc\ at $\xb$, there is $\lambda\in ]0,\eps[$ such that
\begin{equation}\label{sci}
f(\bx)\le \inf f(\xb+\lambda B_X) + \eps.
\end{equation}
Applying Ekeland's variational principle \cite{Eke74}, we find $\xt\in X$ such that
\begin{subequations}
\begin{gather}
\|\xt-\bx\|< \lambda, ~f(\xt)\le f(\bx), \mbox{ and} \label{ekeland1}\\
y\mapsto f(y)+(\eps/\lambda)\|y-\xt\| \mbox{ admits a local minimum at } \xt. \label{ekeland2}
\end{gather}
\end{subequations}
In view of (\ref{ekeland2}), we may apply the Separation Principle at point $\xt$
with the convex Lipschitz function $\varphi:y\mapsto (\eps/\lambda)\|y-\xt\|$ to obtain
a subgradient $\xt^*\in \del f(\xt)$ such that $-\xt^*\in \del \varphi(\xt)$.
From (\ref{ekeland1}) and (\ref{sci}), we derive that $\|\xt-\bx\|\le \eps$ and
$|f(\xt)-f(\xb)|\le\eps$, while from $-\xt^*\in \del \varphi(\xt)$ we get
$\|\xt^*\|\le \eps/\lambda$, so combining with (\ref{ekeland1}) we find
$
\langle \xt^*,\xt-\xb\rangle\le (\eps/\lambda)\lambda=\eps.
$
This completes the proof of the nonemptiness of the sets $\Sd_\eps f(\xb)$.

%%%%%%%%%%%%%%%%%%%%%%%%%%%%%%%%%%%%%%%%%%%%%%%%%%
\ps
We now proceed to the proof of formula (\ref{CDD-formula}).
Let $d\ne 0$ and let $\ld< f'(\xb;d)$.
It suffices to show that
there exists a sequence $((x_n,x_n^*))_n\subset \del f$ verifying
\begin{subequations}
\begin{gather}
x_n\to \xb,\ f(x_n)\to f(\xb),\ \limsup_{n}\,\langle x^*_n,x_n-\xb\rangle\le 0, \mbox{ and} \label{pr1}\\
\ld\le \liminf_{n}\,\langle x^*_n,d\rangle. \label{pr2}
\end{gather}
\end{subequations}

\ps
Let $t_0 \in \,]0,1]$ such that
\begin{equation}\label{etoile}
0\le f(\xb+td)-f(\xb) -\ld t, \quad \forall  t\in [0,t_0]
\end{equation}
and let $z^*\in X^*$ such that 
\begin{equation}\label{HB}
\la z^*,td\ra =-\ld t, \quad \forall t\in\R \quad 
\mbox{and}\quad \|z^*\|={|\ld|}/{\|d\|}.
\end{equation}
Set $K:= [\xb, \xb+t_0d]$ and $g:=f+z^*$. Then, (\ref{etoile}) can be rewritten as
\[
g(\xb) \le g(x), \quad \forall x \in K.
\]
Let $\delta>0$
such that $g$ is bounded from below on $B_{\delta}(K)$.
For each positive integer $n$ such that $1/n<\delta$, let $r_n>0$ such that 
\begin{equation*}\label{penal2}
g(\xb)-1/n^2<\inf_{B_{r_n}(K)}g
\end{equation*}
and let then $\alpha_n>0$ such that
\begin{equation*}\label{penal3}
\inf_{B_{r_n}(K)}g\le \inf_{B_{\delta}(K)}g+\alpha_nr_n.
\end{equation*}
It readily follows from these inequalities that
\begin{equation}\label{penal}
g(\xb) \le g(x) +\alpha_n d_K(x) + 1/n^2, \quad \forall x\in B_{\delta}(K).
\end{equation}
\ps
Applying Ekeland's variational principle to the function
$x\mapsto g(x) +\alpha_n d_K(x)$, we obtain
a sequence $(x_n)_n\subset B_{\delta}(K)$ such that
\begin{subequations}
\begin{gather}
\|\xb-x_n\|<1/n, \label{ek1}\\
g(x_n) + \alpha_n d_K(x_n) \le g(\xb) \label{ek2}\\
x\mapsto g(x) +\alpha_nd_K(x)
+ (1/n) \| x-x_n \|\mbox{ admits a local minimum at } x_n. \label{ek3}
\end{gather}
\end{subequations}
In view of (\ref{ek3}), we may apply the Separation Principle
with the given $f$ and $\varphi= z^*+ \alpha_n d_K+(1/n)\|.-x_n\|$ to obtain 
$x_n^*\in \del f(x_n)$, $\zeta_n^* \in \del d_K(x_n)$ and $\xi^*\in B_{X^*}$ such that 
\begin{equation}\label{decomp}
x_n^*= -z^*-\alpha_n\zeta_n^*-(1/n)\xi^*.
\end{equation}
We show that the sequence $((x_n,x^*_n))_n\subset \del f$
satisfies (\ref{pr1}) and (\ref{pr2}).
\ps
\textit{Proof of (\ref{pr1}).} Combining (\ref{HB}) and (\ref{ek2}) we get
\begin{equation}\label{ek22}
f(x_n)\le f(\xb)+\la z^*,\xb-x_n\ra\le f(\xb)+({|\ld|}/{\|d\|})\|\xb-x_n\|.
\end{equation}
Since $f$ is lower semicontinuous, (\ref{ek1}) and (\ref{ek22}) show
that $x_n\to \xb$ and $f(x_n)\to f(\xb)$.
On the other hand, since $\la\zeta_n^*,x-x_n\ra \le d_K(x)-d_K(x_n)\le 0$ for all $x\in K$, 
it follows from (\ref{HB}) and (\ref{decomp}) that
\begin{equation*}\label{form1}
\la x_n^*,\xb-x_n\ra \ge \la z^*,x_n-\xb\ra - (1/n)\la\xi^*,\xb-x_n\ra
                 \ge -({|\ld|}/{\|d\|})\|x_n-\xb\|-(1/n)\|x_n-\xb\|,
\end{equation*}
showing that $\limsup_n \la x_n^*,x_n-\xb\ra \le 0$.
\ps
\textit{Proof of (\ref{pr2}).} From (\ref{HB}) and (\ref{decomp}) we derive
\begin{equation}\label{form2}
\la x_n^*,d\ra=\la-z^*,d\ra -(1/n)\la \xi^*,d\ra -\alpha_n\la\zeta_n^*,d\ra
\ge\ld-(1/n)\|d\|-\alpha_n\la\zeta_n^*,d\ra.
\end{equation}
We claim that $\la\zeta_n^*,d\ra\le 0$. Indeed, let $P_K{x_n}\in K$
such that $\|x_n -P_K{x_n}\|= d_K(x_n)$.
We have $P_K{x_n}\not= \xb + t_0d $ for large $n$ since $x_n\to \xb$, so there 
exists $t_n>0$ such that $t_0d=\xb+t_0d-\xb=t_n(\xb+t_0d-P_K{x_n})$.
It follows that
\begin{eqnarray*}
({t_0}/{t_n})\la\zeta_n^*,d\ra =\la \zeta_n^*,\xb+t_0 d-P_K{x_n}\ra
&=&\la \zeta_n^*,\xb+t_0d-x_n\ra+ \la\zeta_n^*,x_n-P_K{x_n}\ra\\
&\le&- d_K(x_n) + \|x_n -P_K{x_n}\|=0.
\end{eqnarray*}
Hence $\la\zeta_n^*,d\ra\le 0$. We therefore conclude from (\ref{form2}) that 
$\ld\le \liminf_{n}\,\langle x^*_n,d\rangle$.
This completes the proof.
\cqfd
\ps\noindent
{\em Remarks.} 1. For convex functions, the inequality in (\ref{CDD-formula}) becomes an equality, and the formula is due to Taylor \cite{Tay73} and Borwein \cite{Bor82}.
\ps
2. For elementary subdifferentials, formula (\ref{CDD-formula}) was proved to hold
%(proximal, Fr\'echet, Hadamard, G\^ateaux, \ldots) 
in appropriate Banach spaces, see \cite{JL12}.
The arguments there were based on a specific property of these subdifferentials with
respect to exact inf-convolutions of two functions.
Such an argument is avoided here: formula (\ref{CDD-formula}) is valid in a Banach space $X$
as soon as the subdifferential satisfies the Separation Principle in this space.

%%%%%%%%%%%%%%%%%%%%%%%%%%%%%%%%%%%%%%%%%%%%%%%%%%%%%%%%%%

\section{First-Order Tests for Optimality}

The following lemma comes from our paper \cite{JL12}.
It establishes a mean value inequality using the directional derivative.
For the sake of completeness, we recall the proof.

\begin{lemma}\label{mvi}
Let $X$ be a Hausdorff locally convex space, $f:X\to\xR$ be \lsc, $\xb\in X$ and $x\in\dom f$. Then, for every real number $\lambda\le f(\xb)-f(x)$, there
exist $t_0\in [0,1[$  and $x_0:=x+t_0(\xb-x)\in [x,\xb[$
such that
$\lambda\le f'(x_0;\xb-x)$ and $f(x_0)\le f(x)+t_0\lambda.$
\end{lemma}

\proof
Let $\lambda\in\R$ such that
$\lambda\le f(\xb)-f(x)$ and define $g:[0,1]\to \xR$ by $g(t):=f(x+t(\xb-x))-t \lambda$.
Then $g$ is lower semicontinuous
on the compact $[0,1]$ and $g(0)=f(x)\le f(\xb)-\lambda=g(1)$. Hence $g$ attains its
minimum on $[0,1]$ at a point $t_0\ne 1$. Let $x_0:=x+t_0(\xb-x)\in [x,\xb{[}$.
Then, $f(x_0)-t_0\lambda=g(t_0)\le g(0)=f(x)$ and, since $g(t_0+t)\ge g(t_0)$
for every $t\in {]}0,1-t_0]$, we derive that
\[
\forall t\in {]}0,1-t_0],\quad \frac{f(x_0+t(\xb-x))-f(x_0)}{t} \ge \lambda.
\]
Passing to the lower limit as $t\searrow 0$, we get $f'(x_0;\xb-x)\ge \lambda$.
The proof is complete.
\cqfd

We deduce easily from Lemma \ref{mvi}
a first-order necessary and sufficient condition for
optimality of \lsc\ functions in terms of the directional derivative:

\begin{proposition}\label{mainprop}
Let $X$ be a Banach space, $f:X\to\xR$ be \lsc, $C\subset X$ be convex and $\xb\in C$.
Then, the following are equivalent:\ps
{\rm(a)} $f(\xb)\le f(y)$ for every $y\in C$,\ps
{\rm(b)} $f(\xb)\le f(y)$ for every $y\in C$ such that $f'(y;\xb-y) > 0$.
\end{proposition}
\proof
Obviously, (a) implies (b). We prove that $\neg$(a) implies $\neg$(b).
Assume there exists $x\in C$ such that $f(\xb)> f(x)$.
We must show that there exists $x_0\in C$ such that $f(\xb)> f(x_0)$ and $f'(x_0;\xb-x_0) > 0$.
Applying Lemma \ref{mvi} with $0<\lambda<f(\xb)-f(x)$, we derive that
there exist $x_0\in [x,\xb[\subset C$ and $t_0\in [0,1[$
such that
\[f'(x_0;\xb-x)\ge \lambda> 0 \quad\mbox{and}\quad f(x_0)\le f(x) +t_0\lambda.\]
Since $f(x) +\lambda< f(\xb)$, $\lambda> 0$ and $t_0\in [0,1[$, we have $f(x_0)< f(\xb)$,
and, since $\xb-x_0=t(\xb-x)$ for some $t>0$, we have $f'(x_0;\xb-x_0)=tf'(x_0;\xb-x)> 0$.
\cqfd

The following first-order sufficient condition is a straightforward consequence.
It clearly contains the result quoted in the introduction.
We refer to \cite{CGR04} for a characterization of the solution set of
Minty variational inequalities governed by directional derivatives.
\begin{cor}\label{OCb0}
Let $X$ be a Banach space, $f:X\to\xR$ be \lsc, $C\subset X$ be convex and $\xb\in C$.
Then:
\begin{equation*}\label{sufcond0}
\forall y\in C,\ f'(y;\xb-y) \le 0 \Longrightarrow
\forall y\in C,\ f(\xb)\le f(y).
\end{equation*}
\end{cor}

Now, coming back to our objective, we combine Lemma \ref{mvi} %(or Proposition \ref{mainprop})
with Theorem \ref{link} to establish a first-order necessary and sufficient condition for
optimality of \lsc\ functions in terms of
any subdifferential satisfying the Separation Principle.
This complements our previous result \cite[Theorem 4.2]{JL12} which concerned only
the elementary subdifferentials.

\begin{theorem}\label{main}
Let $X$ be a Banach space, $f:X\to\xR$ be \lsc, $U\subset X$ be open convex and $\xb\in U$.
Then, the following are equivalent:\ps
{\rm(a)} $f(\xb)\le f(y)$ for every $y\in U$,\ps
{\rm(b)} $f(\xb)\le f(y)$ for every $y\in U$ such that $\sup \la \del f(y),\xb-y\ra > 0$.
\end{theorem}
\proof
Obviously, (a) implies (b). Conversely, we show that $\neg$(a) implies $\neg$(b).
We know from Proposition \ref{mainprop} that $\neg$(a) implies the existence of 
$x_0\in U$ such that $f(\xb)> f(x_0)$ and $f'(x_0;\xb-x_0)> 0$.
Let $\eps>0$ such that $x_0+\eps B_X\subset U$, $f(x_0)+\eps< f(\xb)$ and
$f'(x_0;\xb-x_0)> \eps$,
and apply formula (\ref{CDD-formula}) at point $x_0$
and direction $d=\xb-x_0$.
We obtain a pair $(y_\eps,y^*_\eps)\in \del f$ such that
\[
\|y_\eps-x_0\|< \eps, \quad f(y_\eps)<f(x_0)+\eps,
\quad\la y^*_\eps, y_\eps-x_0\ra <\eps,
\quad\la y^*_\eps, \xb -x_0\ra >\eps, 
\]
from which we derive that $y_\eps\in U$, $f(y_\eps)<f(\xb)$ and
$\la y^*_\eps,\xb -y_\eps\ra>0$.
%=\la y^*_\eps,\xb -x_0\ra + \la y^*_\eps,x_0 -y_\eps\ra>\eps-\eps=0$.
\cqfd

As above, the following first-order sufficient condition is a straightforward consequence:

\begin{cor}\label{OCb}
Let $X$ be a Banach space, $f:X\to\xR$ be \lsc, $U\subset X$ be open convex and $\xb\in U$.
Then:
\begin{equation}\label{sufcond}
\forall y\in U,\ \sup \la \del f(y),\xb-y\ra \le 0
\Longrightarrow
\forall y\in U,\ f(\xb)\le f(y).
\end{equation}
\end{cor}

We recall from \cite{JL12} an interpretation of the sufficient condition (\ref{sufcond})
in terms of the monotonic behaviour of the subdifferential operator.
Given a set-valued operator $T:X\rightrightarrows X^*$,
or graph $T\subset X\times X^*$, we let
\[
T^0:=\{(x,x^*)\tq\langle y^*-x^*,y-x\rangle \geq 0~~\forall (y,y^*)\in T\}
\]
be the set of all pairs $(x,x^*)\in X\times X^*$ {\it monotonically related} to $T$.
The operator $T$ is said to be
{\it monotone} provided $T\subset T^0$,
{\it monotone absorbing} provided $T^0\subset T$,
{\it maximal monotone} provided $T= T^0$.

\begin{theorem}\label{mono_absorb}
Let $X$ be a Banach space and let $f:X\to\xR$ be \lsc.
Then, the operator $\del f:X\rightrightarrows X^*$ is monotone absorbing.
In particular, for convex $f$, $\del f$ is maximal monotone.
\end{theorem}
\proof
Let $(x,x^*)\in (\del f)^0$. Then,
$\la y^*-x^*,y-x\ra\ge 0$ for every $(y,y^*)\in\del f$.
Since $\del(f-x^*)(y)=\del f(y) -x^*$,
this can be written as
\begin{equation*}\label{MMg}
\forall (y,y^*)\in\del (f-x^*),\quad \la y^*,y-x\ra\ge 0,
\end{equation*}
so, by Corollary \ref{OCb}, $(f-x^*)(x)\le (f-x^*)(y)$ for every $y\in X$.
We then conclude from the Separation Principle that $x^*\in \del f(x)$.
Thus, $(\del f)^0\subset \del f$.

For convex $f$, $\del f$ is monotone, hence maximal monotone from what precedes.
\cqfd
\end{document}